\documentclass[12pt, final]{amsart}

\usepackage{amsthm, amsmath, amssymb}
\usepackage{subfig}

\usepackage[notref,notcite]{showkeys}
\usepackage[pdftex]{graphicx}
\usepackage{comment}

\newtheorem{theorem}{Theorem}[section]

\newtheorem{lemma}{Lemma}[section]
\newtheorem{proposition}{Proposition}[section]
\newtheorem{definition}{Definition}[section]
\newtheorem{remark}{Remark}[section]
\newtheorem{example}{Example}
\usepackage{color}

\begin{document}

\title[Numerical treatment to a parabolic non-local free boundary problem ]{Numerical treatment to a  non-local  parabolic free boundary problem arising in financial bubbles}


\author{A. Arakelyan}
\address[A.~Arakelyan]{Institute of Mathematics,  NAS of Armenia \\0019 Yerevan, Armenia}
\email[Avetik Arakelyan]{arakelyanavetik@gmail.com}

\author{R. Barkhudaryan}
\address[R.~Barkhudaryan]{Institute of Mathematics, NAS of Armenia\\ 0019 Yerevan, Armenia}
\email[Rafayel Barkhudaryan]{rafayel@instmath.sci.am}

\author{H. Shahgholian}
\address[H.~Shahgholian]{Department of Mathematics\\
	The Royal Institute of Technology\\
	100 44 Stockholm, Sweden}
\email[Henrik Shahgholian]{henriksh@math.kth.se}

\author{M. Salehi}
\address[M.~Salehi]{Department of Mathematics, Statistics and Physics,\\
	Qatar University, P.O. Box 2713, Doha, Qatar}
\email[Mohammad Mahmoud Salehi]{salehi@qu.edu.qa }

\date{}

\begin{abstract}
In this paper we continue to study a non-local free boundary problem arising in financial bubbles. We focus on the parabolic counterpart of the bubble problem and  suggest  an iterative algorithm  which consists of  a sequence of parabolic obstacle problems at each step to be solved, that in turn gives the  next  obstacle function  in the iteration. The convergence of the proposed  algorithm is proved. Moreover, we consider the finite difference scheme for this algorithm and obtain its convergence. At the end of the paper we present  and discuss computational results. \\
%
\end{abstract}

\maketitle

\section{Introduction}

In this paper we shall analyze the time dependent(parabolic)  free boundary problem for a financial bubble problem, from a PDE point of view. Hence the model equation, studied in this paper,  is  the following free boundary problem formulated as a Hamilton-Jacobi equation:
\begin{equation}\label{main}
	\min(Lu,u(t,x)-u(t,-x)-\psi(t,x))=0,\quad (t,x)\in\mathbb{R^+}\times\Omega,
\end{equation}
where $\Omega\subset \mathbb{R}$  is a symmetric bounded domain such that  if $x\in\Omega$ then $-x\in\Omega$  and $\psi\in C^2(\mathbb{R^+}\times\Omega)$. 

 As mentioned above we consider time depended parabolic case, i.e. the operator $L$ is the following parabolic operator
\[
Lu=u_t+Mu
\]
where $M$ is the elliptic counter part of the operator $L$ as defined below:
\[
Mu=a^{ij}(x)D_{ij}u+b^i(x)D_iu+c(x)u,\qquad a^{i,j}=a^{j,i}.
\]
Here the coefficients $a^{i,j}$, $b^i$, $c$ are  assumed to be  continues  and the matrix $[a^{i,j}(x)]$ is positive definite  for all $x\in\Omega$. Additionally we assume that the  coefficients are ``symmetric'' in the domain $\Omega$ i.e. the operator applied to the function $u(-x)$ should be the same as operator applied to the function $u$ at point $-x$:
\begin{equation}\label{sym}
(L\widetilde{u})(x)=(\widetilde{Lu})(x)
\end{equation}
where $\widetilde u(t,x)=u(t,-x)$. The relation \eqref{sym} is the same if we require that $a^{ij}$ and $c$ are even functions  and  $b^i$ is odd.

If the domain $\Omega$ is bounded we are going to consider the problem with the following initial and boundary conditions
\begin{equation}\label{main_bounded}
\begin{cases}
u(0,x)=g(0,x), &x\in\Omega,\\
u(t,x)=g(t,x), &(t,x)\in\mathbb{R^+}\times\partial\Omega,
\end{cases}		
\end{equation}

The specific and important case is the Black-Scholes equation i.e. the domain $\Omega=\mathbb{R}$ and the differential operator is the following
\[M=-\frac{1}{2}\sigma^2 u_{xx}+\rho x u_{x}+ ru.\] Our main concern is to develop numerical method for the following non-local free boundary problem:
\begin{equation}\label{main_bounded-BS}
\begin{cases}
\min(\partial_t u+Mu,u(t,x)-u(t,-x)-\psi(t,x))=0,&(t,x)\in\mathbb{R^+}\times \Omega,\\
u(0,x)=g(0,x), &x\in\Omega,\\
u(t,x)=g(t,x), &(t,x)\in\mathbb{R^+}\times\partial\Omega,
\end{cases}		
\end{equation}
where 
\[\psi(t,x)=x\frac{1 - e^{-(r+\lambda)t}}{r+\lambda}- c\;\;\mbox{and}\;\; \sigma,\rho,r>0.\]

Problem \eqref{main_bounded-BS} arises in modeling of speculative financial bubbles. 
The financial model of speculative trading described in \cite{MR3268054} where it is allowed to profit from other over-valuation and additionally assuming that trading agents may ``agree to disagree''. Due to  speculative trading, asset prices may beat their fundamental values. 

The stationary or finite horizon version of this model was introduced and solved by Scheinkman and Xiong \cite{scheinkman2003}.   Scheinkman and Xiong considers one-dimensional case and they are able to construct an explicit solution based on Kummer functions. 
It was possible to do in one dimensional case as the solution of Black-Scholes equation is possible to express true Kummer function. Other stationary models were studied in \cite{MR3016784, MR2800215}. The multidimensional stationary problem  was considered in \cite{MR3552318} where existence and uniqueness of the viscosity solution were proved.

It is apparent that  that if  $u(t,x)$   is a  solution to equation \eqref{main}, then $u(t,-x)$ is a solution to the reflected problem, with  all ingredients  reflected accordingly.
	 In particular, $u(t,x) \geq u(t,-x) + \psi(t,x) \geq u(t,x) + \psi (t,-x) $, and 
	 $\psi (t,x) + \psi (t,-x) \leq 0 $ is forced as a condition for an  existence theory.
A standing assumption in this paper is that 
the constraint   $\psi$,    and the boundary data $g(t,x)$,  should satisfy the following inequality 
$$\psi(t,x) + \psi(t,-x) \leq 0 ,  \qquad \psi(t,x)\leq g(t,x).$$
For more  details on this problem see \cite{MR3268054}. 

For review of  the obstacle type PDE models in the socio-economic sciences see \cite{MR3268053}, other theoretical aspects of obstacle-type problems you can see in \cite{MR1658612,MR2962060}.

In this paper our  objective is to study, through an increasing iterative algorithm, a solution of the above problem in $[0,T]\times\Omega.$  
The algorithm consists of  a sequence of parabolic obstacle problems at each step that eventually converge to the solution. We also study the corresponding difference scheme developed for the iterative algorithm.

\section{The iterative algorithm}
To deal with the problem  we first  recall the definition of the  so-called viscosity solution following  \cite{MR3268054}. 

\begin{definition}[Viscosity sub/super solution of the equation \eqref{main}]
	
A function $u:\mathbb{R^+}\times\Omega  \to \mathbb{R}$ is a viscosity subsolution (resp. supersolution) of \eqref{main} on  $\mathbb{R^+}\times\Omega$, if $u$ is upper semi-continuous (resp. lower semi-continuous), and if for any function $\varphi\in C^{1,2}(\mathbb{R^+}\times\Omega)$ and any point $(t_0,x_0) \in \mathbb{R^+}\times\Omega$ such that $u(t_0,x_0) = \varphi(t_0,x_0)$ and
\[
u\leq\varphi \text{ (resp. } u \geq \varphi \text{)} \text{ on }  \mathbb{R^+}\times\Omega ,
\]
\[
u\leq g, \text{ (resp. } u\geq g \text{)},\text{ on  the boundary of }  \mathbb{R^+}\times\Omega, 
\]
then
\[
	\min(L\varphi(t_0,x_0), \varphi(t_0,x_0)-\tilde{\varphi}(t_0,x_0)-\psi(t_0,x_0)) \leq 0,
\]
(resp.  $\min(L\varphi(t_0,x_0), \varphi(t_0,x_0) - \tilde{\varphi}(t_0,x_0)-\psi(t_0,x_0)) \geq 0$).
\end{definition}

\begin{definition}[Viscosity solution of the equation \eqref{main}]
A function $u :\mathbb{R}^+\times\Omega\to R$  is a viscosity solution of \eqref{main} on $\Omega$, if and only if $u^*$ is a viscosity subsolution and $u_*$ is a viscosity supersolution on $\Omega$, where
\[
u^*(x)=\limsup_{(l,y)\to (t,x)}u(l,y),
\]
and
\[
u_*(x)=\liminf_{(l,y)\to (t,x)}u(l,y).
\]
\end{definition}

To construct the algorithm, at first let  us define a function $u_0$ (the initial guess) as the solution of the following problem: 
 \[
 	\begin{cases}
 	L u_{0}=0,&(t,x)\in\mathbb{R^+}\times\Omega,\\
	u_{0}(0,x)=	g(0,x), &x\in\Omega,\\
	u_{0}(t,x)=g(t,x), &(t,x)\in\mathbb{R^+}\times\partial\Omega,

 	\end{cases}		
 \]
Inductively, we define the sequence  $\{u_k\}$ by
\begin{equation}\label{iteration}
	\begin{cases}
	\min(L u_{k+1},u_{k+1}-\tilde{u}_{k}-\psi)=0,&(t,x)\in\mathbb{R^+}\times\Omega,\\
	u_{k+1}(0,x)=g(0,x), &x\in\Omega,\\
	u_{k+1}(t,x)=g(t,x), &(t,x)\in\mathbb{R^+}\times\partial\Omega,
	\end{cases}		
\end{equation}
For each $k$ we have an obstacle problem with the obstacle $\tilde{u}_{k}+\psi$. 

In this section our goal is to show that \eqref{iteration} produces a non-decreasing sequence $\{u_k\}$ and then to show that  the sequence $\{u_k\}$  tends to the viscosity solution of \eqref{main}.

\begin{proposition}\label{prop-incr}
The sequence $\{u_k\}$ is non-decreasing.
\end{proposition}
\begin{proof}
It is easy to see that $u_1\geq u_0$, since the function $u_1$ is
the solution of $L u_{1}(t,x)=0$ with obstacle ${u_0+\psi}$,
and $u_0$ is the solution to the same problem without an obstacle.

To prove that $u_2(t,x)\geq u_1(t,x)$ let us examine  equation \eqref{iteration}.
 When $k=1$ the obstacle  is $\Psi_2:=u_1(t,-x)+\psi(t,x)$ and  for the case $k=0$ the obstacle is $\Psi_1:=u_0(t,-x)+\psi(t,x)$.  Since  $u_1(t,x)\geq u_0(t,x)$ we have  $u_1(t,-x)\geq u_0(t,-x)$ and hence  
 $\Psi_2 \geq \Psi_1$ . Furthermore,  $u_2(x)$ and $u_1(x)$ are solutions of the same obstacle problem with obstacles $\Psi_2 $, and $\Psi_1$ respectively. Since $\Psi_2 \geq \Psi_1$ we have by comparison principle (see \cite{MR679313}, page 80, problem 5)   $u_2(t,x)\geq u_1(t,x)$.

By inductive steps we have that $u_{k+1}$ and $u_k$ solve the obstacle problem with obstacle $\Psi_{k+1}:= u_k(t,-x) + \psi$ and $\Psi_k:= u_{k-1}(t,-x) + \psi$ respectively, with $u_{k} \geq u_{k-1}$, and hence $\Psi_{k+1} \geq \Psi_{k}$, and hence by comparison principle $u_{k+1} \geq u_k$.
\end{proof}

We need to prove that the algorithm is bounded above by the solution of \eqref{main}.
\begin{proposition} 
If $w$ is a solution of the problem \eqref{main} then $u_k\leq w$.
\end{proposition} 
\begin{proof}
First of all $u_0\leq w$ since the function $w$ is the solution of $\partial_t w+M w=0$ with obstacle ${\tilde{w}+\psi}$, and $u_0$ is the solution to the same problem without an obstacle.

If we have $u_k(t,x)<w(t,x)$ then by induction we can conclude that $u_{k+1}(t,x)<w(t,x)$ since the function $w$ is the solution of $\partial_t w+M w=0$ with obstacle ${\tilde{w}+\psi}$, and $u_{k+1}$ is the solution to the same problem with smaller obstacle ${\tilde{u}_k+\psi}\leq \tilde{w}+\psi$. 
\end{proof}

\begin{remark}
If the domain $\Omega=\mathbb{R}$ and the operator $L$ is Black-Scholes operator, the existence of solution is proved in \cite{MR3268054} so the algorithm is bounded. If the domain $\Omega$ is bounded, the boundedness of the algorithm is proved in Proposition \ref{prop-bound}.
\end{remark}

\begin{proposition} \label{prop-bound} If the domain $\Omega$ is bounded then the sequence $\{u_k(t,x)\}$ is bounded for every fixed $t$, i.e. there exists a constant $M_t$ such that 
$$
u_k(t,x)\leq M_t,
$$
for all $k\in\mathbb{N}$ and $x\in\Omega$.
\end{proposition}
\begin{proof}

Let $T_0$ be fixed and $h(x)$ be a  symmetric function defined in $\Omega$ and satisfies
$$M h(x) = 1, \quad h(x)\geq 1, \quad \psi(t,x) \leq C_{T_0} h(x)$$
for some large $ C_{T_0}$, such that $\partial_t \psi(t,x)+M\psi(t,x) +  C_{T_0} \geq 1$ holds in $[0,T_0]\times\Omega$; here we have assumed $\psi(t,x) \in C^{1,2}([0,+\infty)\times\Omega)$.
Then from the algorithm defined in \eqref{iteration} we have 
\begin{equation}\label{ineq}
\min\left(\partial_t v_{k+1}(t,x)+Mv_{k+1}(t,x) + 1, v_{k+1}(t,x) - u_{k}(t,-x)-C_{T_0}h(x)\right) \leq 0, 
\end{equation}
where $v_k = u_k - \psi + Ch$. 
Now suppose by induction that 
\begin{equation}\label{bound}
\max (\sup v_k(t,x), \sup (v_k(t,x) + \psi(t,x))) \leq M_{T_0},\;\;\mbox{in}\;\; [0,T_0]\times\Omega,
\end{equation}
 where 
 $$
 M_{T_0}:=\left\{\underset{[0,T_0]\times\Omega}\sup (|g(t,x)|+|\psi(t,x)|+|C_{T_0}h(x)|), \underset{[0,T_0]\times\Omega}\sup v_0(t,x)\right\}.
 $$
  The estimate \eqref{bound} is obviously true for the starting value $u_0$. Let further the maximum value of $v_{k+1}$ be achieved at a point $(t^*,x^*) \in [0,T_0]\times\bar\Omega$.  If it is attained on the boundary then we are done. If it is attained inside the domain, then by the ellipticity of the operator  (and concavity of the graph for $v_{k+1}$ at $(t^*,x^*)$) we have $Mv_{k+1} (t^*,x^*)-\partial_t v_{k+1} (t^*,x^*)\leq 0$, which implies $\partial_t v_{k+1} (t^*,x^*)+Mv_{k+1} (t^*,x^*) + 1 > 0 $, and hence by the inequality \eqref{ineq} we have 
$v_{k+1} (t^*,x^*)  \leq u_k (t^*,-x^*)  + C_{T_0}h(x^*) = v_k(t^*,-x^*)  + \psi (t^*,-x^*)   \leq M_{T_0} $ where we have used $h(x^*) = h(-x^*)$. 

It remains to prove $v_{k+1}(t,x) + \psi(t,x) \leq M_{T_0}$.
We  make a similar argument for $w_{k+1}(t,x) := v_{k+1}(t,-x) + \psi(,t-x)= u_{k+1} (t,-x) + C_{T_0}h(x)$ which satisfies a similar type of equation, with reflected version of the ingredients
\begin{equation*}
\min(\partial_t w_{k+1}(t,x)+Mw_{k+1}(t,x) + 1, w_{k+1}(t,x) - u_{k}(t,x) -C_{T_0}h(x)-\psi (t,-x)) \leq 0.
\end{equation*}
As before let the maximum value of $w_{k+1}$ be achieved at a point  $(t^*,x^*) \in [0,T_0]\times\bar\Omega$. Obviously, if the maximum is on the boundary then we have the desired estimate. Hence we assume the maximum  is attained inside the domain and  by using the same arguments as above we will have the following
\begin{align*}
	w_{k+1}  (t,x)& \leq u_k  (t,x) + \psi (t,-x) + C_{T_0}h(x)\leq\\& \leq  u_k  (t,x) - \psi (t,x)  + C_{T_0}h(x) = v_k(t,x) \leq M_{T_0},
\end{align*}
where we have used  $\psi (t,x) + \psi(t,-x) \leq 0$. Hence we arrive at 
 $$\max (\sup v_{k+1}(t,x), \sup (v_{k+1}(t,x) + \psi(t,x)) \leq M_{T_0}, $$
 in the inductive steps. This completes the proof.
\end{proof}

\subsection{Convergence of the iterative algorithm}
\begin{theorem}\label{iteration-limit}
If $u^k(t,x)$ is the iterative algorithm  given by \eqref{iteration}, and $u=\underset{k\to\infty}{\lim}u^k,$ then $u$ is a unique  continuous viscosity solution of \eqref{main}. 
\end{theorem}
\begin{proof} 
Having a bounded increasing  sequence of continuous functions, the limit function $u$   is lower semi-continuous, i.e.
\begin{equation}\label{sub-sol}
u(t,x)=u_*(t,x):=\liminf_{(l,y)\to (t,x)}u(l,y).
\end{equation}
We also denote by $u^*$ the upper-semi continuous envelop of $u$, i.e.
\[
u^*(t,x):=\limsup_{(l,y)\to (t,x)}u(l,y).
\]
First we show  that the function $u^*$ is a sub-solution to \eqref{main}. For that purpose let us suppose $u^*$ is not a sub-solution. Then there exists $(t_0,x_0)\in\mathbb{R}^+\times\Omega$ and a polynomial  $P$ of degree two satisfying
\[
P\geq u^*, \quad P(t_0,x_0)=u^*(t_0,x_0)
\]
such that
$$
\min\Big\{\partial_t P+M P,P-\tilde{P} - \psi\Big\}>0.
$$
Assume that the first inequality holds. Then
$$
\partial_t P(t_0,x_0)+M P(t_0,x_0)>0
$$
and
$$
P(t_0,x_0)>\tilde{P}(t_0,x_0)+\psi(t_0,x_0).
$$
Substituting the values for $P(t_0,x_0)$, the last inequality can be rewritten in the following way:
$$
u^*(t_0,x_0)>\tilde{u}^*(t_0,x_0)+\psi(t_0,x_0).
$$
Using continuity of $f$ and $\psi$, and the fact $u_j\uparrow u^*$, we can deduce
that there exists a number $r>0$ such that
\[
\partial_t P(t,x)+M P(t,x) >0,\quad (t,x)\in B((t_0,x_0),r),
\]
\[
u^*(t,x)>\tilde{u}^*(t,x)+\psi(t,x),\quad (t,x)\in B((t_0,x_0),r).
\]
and (using continuity of $u^*-\tilde{u}^*$) there exists a positive number $\mu <u^*(t_0,x_0)-\tilde{u}^*(t_0,x_0)-\psi(t_0,x_0)$ such that
\begin{equation}\label{eq-00}
u_j(t,x)> \tilde{u}_{j-1}(t,x)+\psi(t,x)+\mu,\quad (t,x)\in \times B((t_0,x_0),r).
\end{equation}
Denote
$$
P_\varepsilon=P+\varepsilon\eta(x),
$$
where $\eta(t,x)$ is a function which satisfy $\partial_t\eta- L\eta=-1$, $\eta\geq0$ and $\eta(t_0, x_0)=0$.
If $\varepsilon>0$ is small enough, then
\begin{equation}\label{eq-0}
\partial_t P_\varepsilon(t,x)+M P_\varepsilon(t,x)>0,\quad (t,x)\in \times B((t_0,x_0),r).
\end{equation}
Next observe that
$$
P_\varepsilon(t,x)> P(t,x) \ge u^*(t,x),\quad (t,x)\neq (t_0,x_0).
$$
As $u_k\uparrow u^*$, we can choose $j$ large enough to satisfy
\begin{equation}\label{eq-1}
\inf_{B((t_0,x_0),r)}\left(P_\varepsilon-u_j\right)<\min_{\partial B((t_0,x_0),r)}\left(P_{\varepsilon}-u^*\right)
\end{equation}
and
\begin{equation}\label{eq-2}
\inf_{B((t_0,x_0),r)}\left(P_\varepsilon-u_j\right)<\mu.
\end{equation}
Take $Q_\varepsilon=P_\varepsilon-c$, where $c$ is a constant chosen in such a way, that $Q_\varepsilon$ touches $u_j$ from above at some $(t',x')\in B((t_0,x_0),r)$ (the inequality \eqref{eq-1} guarantees that the first touch point in $B((t_0,x_0),r)$ will be not on the boundary of $B((t_0,x_0),r)$).

We have constructed at this point a function $Q_\varepsilon$ satisfying
 the following conditions:
\begin{equation}\label{eq-3}
Q_\varepsilon (t,x)\geq u_j(t,x),\quad (t,x)\in B((t_0,x_0),r)
\end{equation}
\begin{equation}\label{eq-4}
Q_\varepsilon (t',x')=u_j(t',x'),\quad where\quad (t',x')\in B((t_0,x_0),r).
\end{equation}

Since $u_j$ is a viscosity subsolution (and, in fact, a solution) of
$$
\min\left\{\partial_t u_j +Mu_j,u_j-\tilde{u}_{j-1}-\psi\right\}=0,\quad u_j|_{\partial \Omega}=g_1,
$$
then, by the definition of viscosity subsolution and \eqref{eq-3}-\eqref{eq-4}, we obtain
\begin{equation}\label{eq-5}
\min\left\{\partial_tQ_\varepsilon(t',x')+M Q_\varepsilon(t',x'), Q_\varepsilon(t',x')-\tilde{u}_{j-1}(t',x')-\psi(t',x')\right\}\le0 .
\end{equation}
Using \eqref{eq-0}, we have $\partial_t Q_\varepsilon (t',x')+M Q_\varepsilon (t',x')>0$, so the only possibility to satisfy \eqref{eq-5} is
$$
Q_\varepsilon(t',x')\le \tilde{u}_{j-1}(t',x')+\psi(t',x').
$$
This means that
$$
P(t',x')\le \tilde{u}_{j-1}(t',x')+\psi(t',x') +c - \varepsilon \eta(t',x').
$$
Then
$$
u_j(t',x')\le u^*(t',x')\le P(t',x')\le \tilde{u}_{j-1}(t',x')+\psi(t',x') +c - \varepsilon \eta(t',x'),
$$
hence
$$
u_j(t',x')\le \tilde{u}_{j-1}(t',x')+\psi(t',x') +c.
$$
But, we deduce from \eqref{eq-00}
$$
u_j(t',x')> u_{j-1}(t',x')+\psi(t',x') +\mu.
$$
This is a contradiction, since by \eqref{eq-2}, it follows $c<\mu$. Hence  $u^*$ is a sub-solution of \eqref{main}.

Let us now discuss the super-solution properties of  $u=u_*$, see \eqref{sub-sol}. Suppose  $u$ is not a super-solution. Then there exists $(t_0,x_0)\in\mathbb{R}^+\times\Omega$ and a polynomial $P$ of degree two satisfying
$$
P\leq u, \quad P(t_0,x_0)=u(t_0,x_0)
$$
such that
\[
\min\Big\{\partial_t P(t_0,x_0)+M P(t_0,x_0),P(t_0,x_0)-\tilde{P}(t_0,x_0) - \psi(t_0,x_0)\Big\}<0.
\]
 Then
\[
\partial_t P(t_0,x_0)+M P(t_0,x_0)<0
\]
or
\begin{equation}\label{super-ineq-P}
P(t_0,x_0)<\tilde{P}(t_0,x_0)+\psi(t_0,x_0).
\end{equation}
Let us consider the first inequality. Using continuity of $u_j, f$ we can deduce that there exists a number $r>0$  such that
\[
\partial_t P(t,x) +M P(t,x)<0,\quad (t,x)\in B((t_0,x_0),r).
\]
Like in the previous case we will construct new polynomial $Q=P-c$ which will touch $u_j$ at some point $x'\in B((t_0,x_0),r)$ i.e.
\begin{equation*}
Q_\varepsilon (t,x)\leq u_j(t,x),\quad (t,x)\in B((t_0,x_0),r)
\end{equation*}
\begin{equation*}
Q (t',x')=u_j(t',x'),\quad \text{where}\quad (t',x')\in B((t_0,x_0),r).
\end{equation*}

Since $u_j$ is a viscosity supersolution (and, in fact, a solution) of
\begin{equation}\label{super-uj-obst}
\min\left\{\partial_t u_j +Mu_j, u_j- \tilde{u}_{j-1}-\psi\right\}=0,\quad u_j|_{\partial \Omega}=g_1,
\end{equation}
we will get contradiction as $\partial_t Q+MQ<0$.

It remains to show that inequality \eqref{super-ineq-P} also cannot be hold. For that purpose let us substitute the values for $P$ and rewrite \eqref{super-ineq-P}:
\[
u(t_0,x_0)-\tilde{u}(t_0,x_0)<\psi(t_0,x_0).
\]
The function $u-\tilde{u}$ is continuous and if the value of $j$ is enough big then 
\[
u_j(t,x)-\tilde{u}_{j-1}(t,x)<\psi(t,x),\quad (t,x)\in B((t_0,x_0),r).
\]
This is a contradiction as $u_j$ should satisfy \eqref{super-uj-obst}

The continuity of $u$ follows from the comparison principle (see \cite{MR3268054}). Indeed, it follows from the comparison principle  that the super solution should be greater or equal to the subsolution, but from the definition of $w$ it follow that $u^*\geq u$, so $u=u^*$ is a continuous viscosity solution of \eqref{main}.
\end{proof}

\section{Finite difference scheme for the iterative algorithm }
For every step of the above algorithm we should solve an obstacle problem and we are going to use finite difference scheme  to do this numerically. The finite difference scheme was extensively used for numerical solutions of variational inequalities, one-phase obstacle problems of elliptic and parabolic type, and in particular, for valuation of American type option (for details, see \cite{WDH} and references in these papers). 

In 2009, the explicit finite difference scheme has been applied for one-dimen\-sio\-nal parabolic obstacle problem in connection with valuation of American type options (see \cite{MR2532350}). It has been proved, that under some natural conditions, the finite difference scheme converges to the exact solution and the rate of convergence is $O(\sqrt{\Delta t}+\Delta x)$. Here $\Delta x$ and $\Delta t$ are space- and time- discretization steps. Recently in the works  \cite{MR2961456, MR2894333Trans, ajm689, MR222222} finite difference scheme and the convergence results have been applied for the one-phase and two-phase elliptic obstacle problems.

In this section we assume that our bubble problem \eqref{main_bounded-BS} is defined on $\Omega^T=[0,T]\times \Omega,$ where $\Omega=[-a,a],$ and $M$ is taken the Black-Scholes operator as defined in introduction.

To construct a finite difference scheme we start by discretizing the domain $\Omega^T=[0,T]\times\Omega$ into a regular uniform mesh. We will denote by $\Omega_h$ and $\Omega_h^T$ the uniform discretized sets of    $\Omega$ and $\Omega^T$ respectively.  For the sake of convenience we set  $h$ as a shorthand of a pair $(\Delta x, \Delta t)$.   
Thus, 
$$
\Omega_h^T=\{(n\Delta t,-a+ {m}\Delta x )\in \mathbb{R}^2,\;\; n=0,1,2,\dots,N\;\;\mbox{and}\;\; m=0,1,2,\dots,M\},
$$ 
where $\Delta t=T/N,$ and $\Delta x=2a/M$.

The discrete Black-Scholes operator is defined as follows 
\begin{multline}\label{disc-black-scholes}
	L_hv(t,x)\equiv\frac{v(t,x)-v(t-\Delta t,x)}{\Delta t}-\frac{\sigma^2}{2}\cdot\frac{v(t,x+\Delta x)-2v(t,x)+v(t,x-\Delta x)}{(\Delta x)^2}+\\\\+ \rho x \frac{v(t,x+\Delta x)-v(t,x-\Delta x)}{2\Delta x}+rv(t,x),
\end{multline}
for any  interior point $(t,x)\in\Omega_h^T$. 

   Let $u^{(k)}=u^{(k)}(t,x)$ be a solution to the iterated obstacle problem with obstacle $u^{(k-1)}(t,-x)+\psi(t,x)$. By  $u^{(k)}_h$  we set the solution  to the following nonlinear system: 
\begin{equation}\label{iteration-discrete}
\begin{cases}
\min\left(L_h u_h^{(k)}(t,x),u_h^{(k)}(t,x)-u^{(k-1)}(t,-x)-\psi(t,x)\right)=0, &(t,x)\in\Omega_h^T,\\
u_h^{(k)}(0,x)=g(0,x), &x\in\Omega_h,\\
u_h^{(k)}(n\Delta t,\pm a)=g(n\Delta t,\pm a), & n=0,1,2,\dots,N.\\
\end{cases}		
\end{equation}

We set  the variational form of the parabolic obstacle problem
\[
F_h^{g}[v]\equiv\min\left(L_h v(t,x),v(t,x)-g(t,x)\right).
\]
Then the following discrete comparison principle for the difference schemes holds.
\begin{lemma}\label{comp-principle}
	Let $L_h$ defined by \eqref{disc-black-scholes} satisfying $|\rho x|\le \frac{\sigma^2}{\Delta x},$ for every $x\in \Omega,$ where $\Omega=[-a,a]$. If $u(t,x)$ and $v(t,x)$ are piecewise continuous functions and satisfy 
$$
F_h^{g}[u]\ge F_h^{g}[v], \;\;\mbox{for all}\;\; (t,x)\in [0,T]\times \Omega,
$$
$$
 u(0,x) \ge v(0,x), \;\;\mbox{for all}\;\; x\in \Omega_h,
$$
then 
$$
u(t,x) \ge v(t,x), \;\;\mbox{for all}\;\; (t,x)\in \Omega_h^T.
$$
\end{lemma}
\begin{proof}

We shall prove by induction
\begin{equation}\label{comparison-induction}
u(n\Delta t,x) \ge v(n\Delta t,x), \;\;\forall x\in\Omega_h,
\end{equation}
for all $n=0,1,2,\dots,N,$ where $N=T/\Delta t$. In the case $n=0$ the inequality \eqref{comparison-induction} coincides with the lemma assumption.  Assume that \eqref{comparison-induction} holds for $n=k,$ we shall prove that it holds for $n=k+1$ as well. We set $t^{k+1}\equiv(k+1)\Delta t$.  For $(t^{k+1},x)$ with $x\in \Omega_h,$ if $v(t^{k+1},x)-g(t^{k+1},x)=F_h^{g}[v](t^{k+1},x),$ then clearly  we get
$$
u(t^{k+1},x)-g(t^{k+1},x)\geq F_h^{g}[u](t^{k+1},x)\geq F_h^{g}[v](t^{k+1},x)=v(t^{k+1},x)-g(t^{k+1},x),
$$
which implies $u(t^{k+1},x)\geq v(t^{k+1},x)$ in this case.	 On the other hand, if $L_h v(t^{k+1},x)=F_h^{g}[v](t^{k+1},x),$ then 
\begin{equation}\label{BS-scheme-ineq}
L_h u(t^{k+1},x)\geq F_h^{g}[u](t^{k+1},x)\geq F_h^{g}[v](t^{k+1},x)=L_h v(t^{k+1},x),
\end{equation}
for every $x\in\Omega_h$. In the sequel we use the following notation: $$w_m^{k}\equiv w(k\Delta t,-a+m\Delta x)\equiv w(t^k,x_m),$$
for all $k=1,2,\dots,N$ and $m=0,1,2,\dots,M$.
In view of \eqref{BS-scheme-ineq} we have
$$
0\leq L_h[u(t^{k+1},x)-v(t^{k+1},x)]=L_h[u_m^{k+1}-v_m^{k+1}],
$$
where $x=-a+m\Delta x \in\Omega_h$. Using the definition \eqref{disc-black-scholes} after simple computation one gets
\begin{multline*}
	0\le \Delta t\cdot L_h[u_m^{k+1}-v_m^{k+1}]=e(u_m^{k+1}-v_m^{k+1})+                       \\
	+d_m(u_{m+1}^{k+1}-v_{m+1}^{k+1})+f_m(u_{m-1}^{k+1}-v_{m-1}^{k+1})-(u_{m}^{k}-v_{m}^{k}),
\end{multline*}
where
\[
e=\left(1+\sigma^2\frac{\Delta t}{(\Delta x)^2}+r\Delta t\right),\;\;
d_m=\left(\frac{\rho x_m }{2}\frac{\Delta t}{\Delta x}-\sigma^2\frac{\Delta t}{2(\Delta x)^2}\right),
\]
\[
f_m=\left(-\frac{\rho x_m }{2}\frac{\Delta t}{\Delta x}-\sigma^2\frac{\Delta t}{2(\Delta x)^2}\right).
\]
Let us rewrite in the matrix form the above equation for all $m=0,1,2,\dots,M$. We have
\begin{equation}\label{implicit-scheme-matrix}
\Delta t\cdot L_h(U^{k+1}-V^{k+1})=A\cdot(U^{k+1}-V^{k+1})-(U^{k}-V^{k}), 
\end{equation}
where $L_h(U^{k+1}-V^{k+1}),U^{k+1}$ and $V^{k+1}$ are column  matrices of $L_h[u_m^{k+1}-v_m^{k+1}],u_m^{k+1}$ and $v_m^{k+1}$ respectively. The matrix $A$ will be  a tridiagonal matrix, such that $A=e I-B$, where $I$ is an identity matrix, $B$ is a tridiagonal matrix with $0$ on the main diagonal and $-d_m$, $-f_m$  on the first diagonals above and below to the main diagonal. According  to \cite[Chapter $6$]{nonnegative-matrices-book}  the matrix $A$ satisfies the properties of an $M$-matrix. Thus, there exists an inverse matrix $A^{-1}$ with non-negative  elements, provided $e>\rho(B),$ where $\rho(B)$ is the spectral radius of the matrix $B$. Let us verify the condition $e>\rho(B)$. To this end,
we observe that $\rho(B)\leq ||B||,$ where the norm  $||.||$ is taken with respect to the rows, i.e. $||B||:=\underset{i}{max}\;\left(\underset{l}{\sum}|b_{i,l}|\right)$. On the other hand, according to the definition of $B$ we get
$$
||B||=\underset{i}{max}\;\left(\underset{l}{\sum}|b_{i,l}|\right)=\underset{i}{max}\;\left(|f_i|+|d_i|\right)=\sigma^2\frac{\Delta t}{(\Delta x)^2}< e,
$$
due to the lemma condition $|\rho x|\le \frac{\sigma^2}{\Delta x}$. Hence,  
$\rho(B)\leq ||B||<e$.

Now,  multiplying by  $A^{-1}$ both sides of the equation \eqref{implicit-scheme-matrix} we arrive at:
\begin{equation}
(U^{k+1}-V^{k+1})=A^{-1}\cdot(U^{k}-V^{k})+\Delta tA^{-1}\cdot L_h(U^{k+1}-V^{k+1}). 
\end{equation}
Recalling that the elements of $A^{-1}$ and $L_h(U^{k+1}-V^{k+1})$ are non-negative we conclude that $u_m^k\geq v_m^k$ implies $u_m^{k+1}\geq v_m^{k+1}$. This completes the proof. 

\end{proof}

\begin{lemma}\label{incr-scheme}
	Let $u_h^{(k)}(t,x)$ be a solution to \eqref{iteration-discrete}, when $|\rho x|\le \frac{\sigma^2}{\Delta x}$ for every $x\in\Omega.$ Then for every $k\in\mathbb{N}$ we have 
	$$ 
	u_h^{(k)}(t,x)\leq u_h^{(k+1)}(t,x),\;\;\mbox{for all}\;\; (t,x)\in\Omega_h^T.
	$$ 
	Moreover, this sequence is  bounded above, which in turn implies its convergence when $k\to\infty$.
\end{lemma}
\begin{proof}
	To prove the statement we apply  the discrete comparison principle for the variational form of the obstacle problems.  In view of Proposition \ref{prop-incr}  we have $u^{(k-1)}(t,-x)\leq u^{(k)}(t,-x)$. 
	This implies
	\begin{multline*}
	0=\min\left\{L_hu^{(k+1)}_h(t,x),\; u^{(k+1)}_h(t,x)-u^{(k)}(t,-x)-\psi(t,x)\right\}=\\
	=\min\left\{L_hu^{(k)}_h(t,x),\; u^{(k)}_h(t,x)-u^{(k-1)}(t,-x)-\psi(t,x)\right\}\geq\\
	\geq
	\min\left\{L_hu^{(k)}_h(t,x),\; u^{(k)}_h(t,x)-u^{(k)}(t,-x)-\psi(t,x)\right\}.
	\end{multline*}
	Thus, one can apply the discrete comparison principle (see Lemma \ref{comp-principle})  for a parabolic obstacle problem with obstacle $u^{(k)}(t,-x)+\psi(t,x)$. This yields
	$u_h^{(k+1)}(t,x)\geq u_h^{(k)}(t,x)$ for all $(t,x)\in\Omega_h^T$.
	
	Let us  prove the boundedness of the sequence $u_h^{(k)}(t,x)$. To do this, we set by 
	$$
	M=\left(\underset{[0,T]\times\overline{\Omega}}\sup |\psi(t,x)|+\underset{[0,T]\times\overline{\Omega}}\sup|w(t,x)|\right)<+\infty,
	$$
	where $\psi(t,x)$ is the obstacle and $w(t,x)$ is a continuous viscosity solution of our parabolic bubble problem \eqref{main_bounded-BS}.
    Then we obtain
    \begin{multline*}
    \min\left\{L_hu^{(k)}_h(t,x),\; u^{(k)}_h(t,x)-u^{(k-1)}(t,-x)-\psi(t,x)\right\}=0\leq\\
    \le
    \min\left\{L_h (M),\; M-u^{(k-1)}(t,-x)-\psi(t,x)\right\}.
    \end{multline*}
    Since $u^{(k)}_h(0,x)\leq M$ for every $x\in\Omega_h,$ then the discrete comparison principle (Lemma \ref{comp-principle}) implies $u^{(k)}_h(t,x)\leq M$ for every $(t,x)\in\Omega_h^T$.

\end{proof}

We are in a position to prove the convergence of the difference scheme for the iterative algorithm.

\begin{proposition}
	Let  $|\rho x|\le \frac{\sigma^2}{\Delta x}$ for every $x\in\Omega,$ and  $w(t,x)$ be a continuous viscosity solution to the parabolic bubble problem \eqref{main_bounded-BS} determined on $[0,T]\times\Omega.$ Define  $u^{(k)}(t,x)$ to be an increasing iterative sequence \eqref{iteration}. If we set by 
    $u^{(k)}_h(t,x)$ the corresponding difference scheme, then 
   $$ \underset{h\to 0}{\lim}\left(\underset{k\to \infty}{\lim}u_h^{(k)}(t,x)\right)=w(t,x).$$
\end{proposition}
\begin{proof}
Assume that $w(t,x)$ is a continuous viscosity solution to the parabolic bubble problem. Due to Lemma \ref{incr-scheme} we know that $u^{(k)}_h(t,x)$ is convergent and therefore there exist some $v_h(t,x)$ such that $u^{(k)}_h\nearrow v_h$ as $k\to \infty$. By Theorem \ref{iteration-limit} it is clear that $v_h$ is a solution to 
\begin{equation}
\begin{cases}
\min\left\{L_h v_h(t,x),\; v_h(t,x)-w(t,-x)-\psi(t,x)\right\}=0, &(t,x)\in\Omega_h^T,\\
v_h(0,x)=g(0,x), &x\in\Omega_h,\\
v_h(n\Delta t,\pm a)=g(n\Delta t,\pm a), & n=0,1,2,\dots,N.
\end{cases}		
\end{equation}

On the other hand $v_h(t,x)$ is a difference scheme for the following parabolic obstacle problem:
\begin{equation}\label{obstacle-problem-v}
\begin{cases}
\min\left\{L v(t,x),\; v(t,x)-w(t,-x)-\psi(t,x)\right\}=0, &(t,x)\in\Omega^T,\\
v(0,x)=g(0,x), &x\in\Omega,\\
v(t,x)=g(t,x), &(t,x)\in [0,T]\times\partial\Omega.
\end{cases}		
\end{equation}
The solution to the above obstacle problem \eqref{obstacle-problem-v} is unique. But $w(t,x)$  also solves \eqref{obstacle-problem-v}, which implies $v(t,x)=w(t,x)$. Now applying Barles-Souganidis theorem  for difference schemes (see \cite{barles-souganidis91})  we obtain $v_h(t,x)\to v(t,x)$ uniformly as $h\to 0$.  This completes the proof.

\end{proof}

Next, we want to estimate  $|w(t,x)-v_h(t,x)|$.  This is nothing else but the difference between the exact solution and the difference scheme for a parabolic obstacle problem \eqref{obstacle-problem-v}. In recent years  there has been given much attention to these type of estimates (see \cite{MR2532350,MR2894333Trans,MR1466804,error-onephase-chinese}). We will mainly follow the above mentioned work \cite{MR2532350}, which considers the problem for  American option valuation. It is worthwhile to mention that they obtain the convergence rate of the order  $O(\Delta x + \sqrt{\Delta t})$.

\begin{proposition}
Let $w(t,x)$ be a  viscosity solution to the parabolic bubble problem \eqref{main_bounded-BS}. If $w(t,x)\in C^{2,1}_{x,t}(\Omega^T),$ then 
$$|w(t,x)-v_h(t,x)|\leq C_\Omega\cdot (\Delta x + \sqrt{\Delta t}),$$
where $v_h(t,x)=\underset{k\to \infty}{\lim}u_h^{(k)}(t,x)$.
\end{proposition}
\begin{proof}
For the proof we recall the parabolic obstacle problem \eqref{obstacle-problem-v}. As we have seen its solution is $w(t,x),$ hence we consider the error analysis for the equation \eqref{obstacle-problem-v} with obstacle $g(t,x)=w(t,-x)+\psi(t,x)$. 
Then we proceed as in \cite{MR2532350}.   
\end{proof}


\section{Numerical results}
In this section we present computational test for the non-local parabolic bubble problem. 
\begin{example}\label{firstexample}
We consider numerical solution of the parabolic financial bubble problem in the domain $\Omega^T=[0,3]\times [-2,2]$.
\begin{equation}\label{main_num}
\min(\partial_t u -(1/2)\sigma^2 u'' + \rho x u' + r u, u-\tilde{u}-\psi)=0,\quad (t,x)\in \Omega^T,
\end{equation}
where $r = 10$, $\rho = 5$, $\sigma = 1$, $\lambda = 1$ and the obstacle function is 
\[\psi(x)=x\frac{1 -e^{-(r + \lambda)*t}}{(r + \lambda)} - 0.001.
\]

\begin{figure}[h]
		\begin{center}
			\subfloat[$u_h$]{\includegraphics[scale=.7]{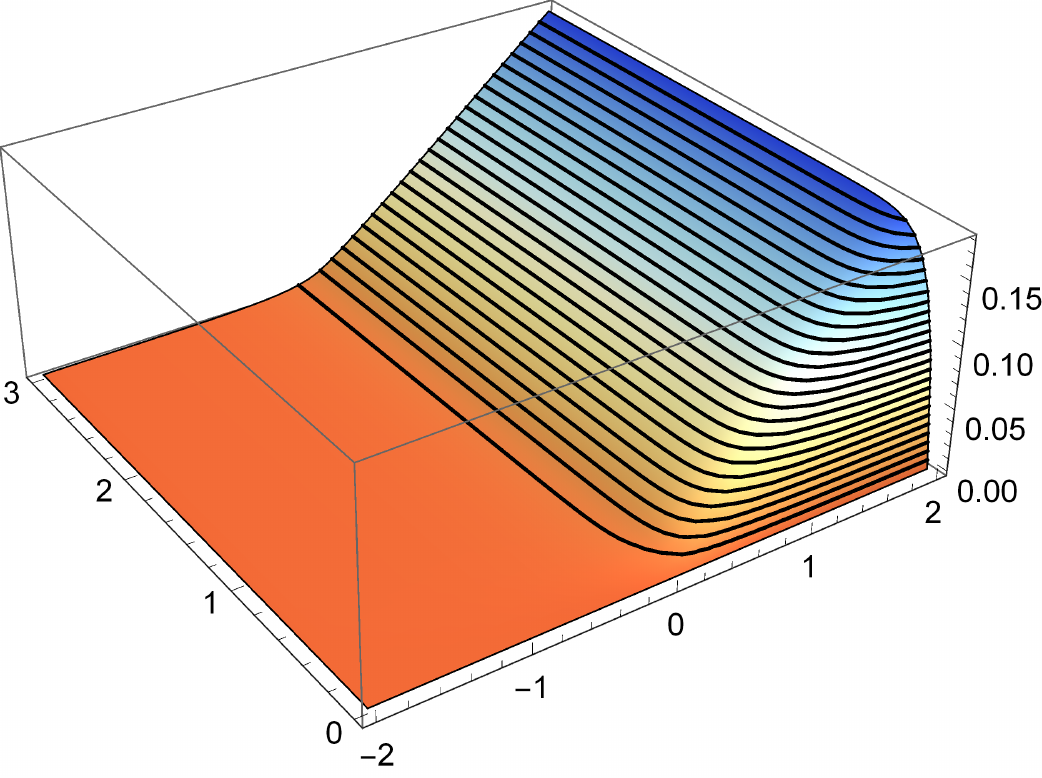}\label{fig1A}} 
			\hspace{.1cm}
			\subfloat[$u_h$]{\includegraphics[scale=.5]{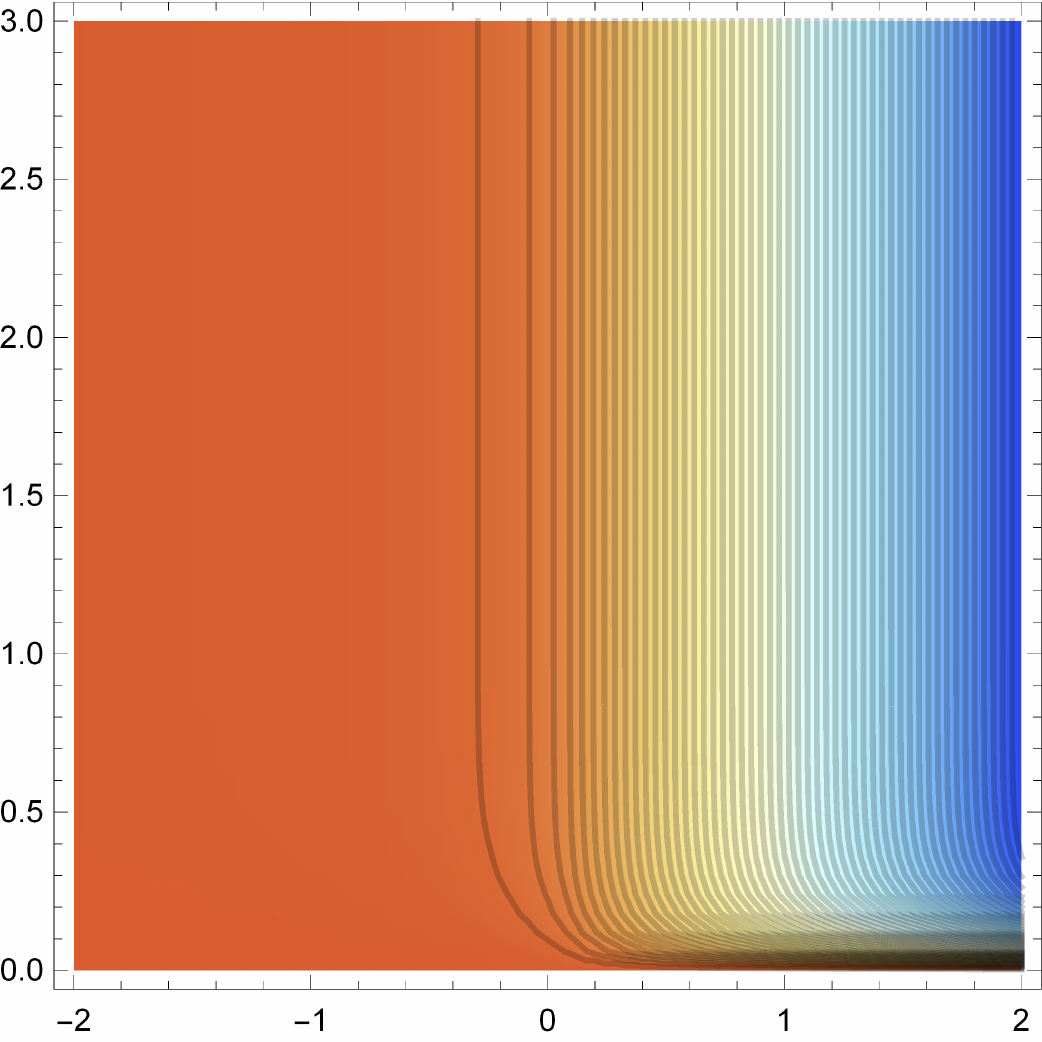}\label{fig1B}}
		\end{center}
		\caption{The numerical solution and the solution densities.}\label{parabolic_num}
\end{figure}
 The numerical solution and its level sets are shown in Figure \ref{parabolic_num} with the use of $50\times 50$ discretization points and after $5$ iterations steps. 

\end{example}

\subsection{One dimensional stationary case}
One dimensional stationary case of financial bubble was considered in \cite{scheinkman2003}, where exact solution was constructed. Following \cite{scheinkman2003} we are going to consider one dimensional stationary case.

 Let 
 \begin{equation}
h(x)=\begin{cases}
U(\frac{r}{2\rho},\frac{1}{2},\frac{\rho}{\sigma^2}x^2)& x\leq 0,\\
\frac{2\pi}{ \Gamma(\frac{1}{2}+(r/2\rho))\Gamma(\frac{1}{2})}M(\frac{r}{2\rho},\frac{1}{2},\frac{\rho}{\sigma^2}x^2)-U(\frac{r}{2\rho},\frac{1}{2},\frac{\rho}{\sigma^2}x^2)& x>0,\\
\end{cases}
 \end{equation}
 where $\Gamma$ is the gamma function, and $U:\mathbb{R}^3\to\mathbb{R}$ is a confluent hypergeometric function of the first kind, $M:\mathbb{R}^3\to\mathbb{R}$ is a confluent hypergeometric function of the second kind. The function $h(x)$
 is positive and increasing in $(-\infty,0)$.
 
Using $h$ function, the exact solution of the bubble problem can be written as
 \begin{equation}\label{exact}
q(x)=\begin{cases}
\frac{b}{h(-k^*)}h(x), & x<k^*,\\
\frac{x}{r+\lambda}+\frac{b}{h(-k^*)}h(-x)-c, & x<k^*,
\end{cases}
 \end{equation}
 where
 \[
 b=\frac{1}{r+\lambda}\frac{h(-k^*)}{h'(k^*)+h'(-k^*)}
 \]
and $k^*$ is a free boundary of the problem which satisfies
\begin{equation} \label{fboundary}
[k^*-c(r+\lambda)][h'(k^*) + h' (-k^*)]-h(k^* ) + h(-k^*) = 0.
\end{equation}
Equation \eqref{fboundary} can be rewritten in simpler form.
\begin{multline}\label{boundarypoint}
U\left(\frac{1}{2} \left(\frac{r}{\rho }-1\right),\frac{1}{2},\frac{{k^*}^2 \rho }{\sigma ^2}\right)\times\\
\times
 \left(2 {k^*}^2 \rho ^2
(k^*-c (\lambda +r))+\sigma^2 (c (\lambda +r) (\rho -r)+k^* (r-2 \rho ))\right)+\\+
k^*(r-2 \rho)
U\left(\frac{1}{2} \left(\frac{r}{\rho }-1\right),\frac{3}{2},\frac{{k^*}^2 \rho }{\sigma ^2}\right) \left(2 k^* \rho 
(k^*-c (\lambda +r))-\sigma ^2\right)=0.
\end{multline}

Next example is devoted to the stationary case of the problem \eqref{main_num}, and we are going to compare exact solution with the numerical solution of the iterative algorithm.
\begin{example}

Here we will consider the finite horizone (i.e. stationary case) for the problem given in Example \ref{firstexample}. 
\begin{equation}
\min( -(1/2) \sigma^2 u'' + \rho x u' + r u,u-\tilde{u}-\psi)=0,\quad x\in\mathbb{R},
\end{equation}
where
\[\psi(x)=\frac{x}{(r + \lambda)} - 0.001.
\]
and $r$, $\rho$, $\sigma$ and $\lambda$  are defined in the previous Example.
\begin{figure}[h]
\begin{center}
	\includegraphics[scale=.6]{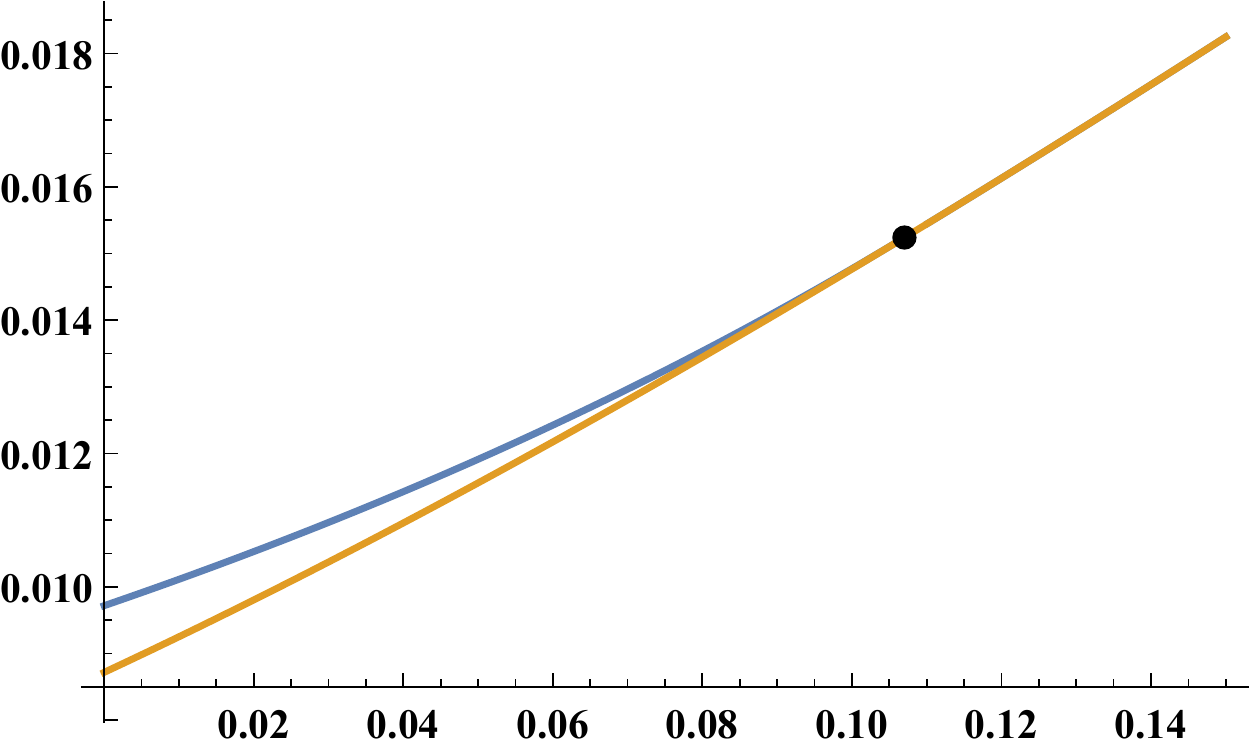}
	\caption{The orange is a solution  $u$, the blue is an obstacle $\tilde{u}+\psi$.}\label{Fig-Bubble-exact}
\end{center}
\end{figure}

Using \eqref{exact} and \eqref{boundarypoint}, the exact solution can be written  as
\[
u(x)=\frac{e^{5 x^2-5 {k^*}^2}}{\sqrt{5 \pi } \left(440 {k^*}^2+44\right)}
\begin{cases}
E_{\frac{3}{2}}\left(5 x^2\right)  & x\leq 0, \\
2\left(\sqrt{5 \pi } x \left(\text{erf}\left(\sqrt{5} x\right)+1\right)+e^{-5 x^2} \right) & x<0 \,\, \&\& \,\, x  \leq k^*,\\
E_{\frac{3}{2}}\left(5 x^2\right)  +\frac{\sqrt{5 \pi } \left(440 {k^*}^2+44\right)}{e^{5 x^2-5 {k^*}^2}} \left(\frac{x}{11}-\frac{1}{10}\right)& x>k^*,
\end{cases}
\]
where
\[
E_{n}\left(x\right)=\int _1^{\infty }\frac{e^{-tx}}{t^n}dt,
\quad
\text{erf}(x)=\frac{2}{\sqrt{\pi }}\int _0^x e^{-t^2}d t
\]
and the point $k^*=0.107028$ is the free boundary point and is the unique real root of the following polynomial 
\[10000 {k^*}^3 - 110 {k^*}^2 -11=0.
\]
In  Figure \ref{Fig-Bubble-exact} the exact solution $u(x)$,  the obstacle function $\tilde{u}+\psi$ and the free boundary point $k^*$ are presented.
\begin{figure}[h]
	\begin{center}
		\subfloat[Error between exact and numerical solution at times $T=0.5$ (orange), $T=1$ (blue dash).]{\includegraphics[scale=.6]{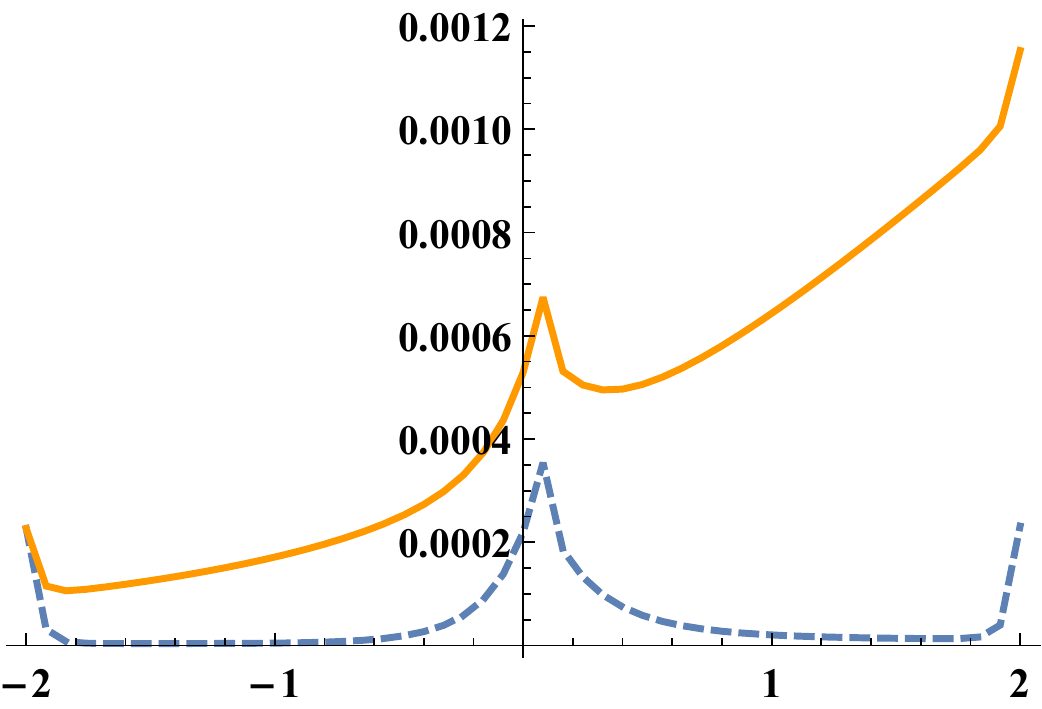}\label{fig2A}} 
		\hspace{.5cm}
		\subfloat[Error between exact and numerical solution at times $T=1$ (blue dash), $T=3$ (blue).]{\includegraphics[scale=.6]{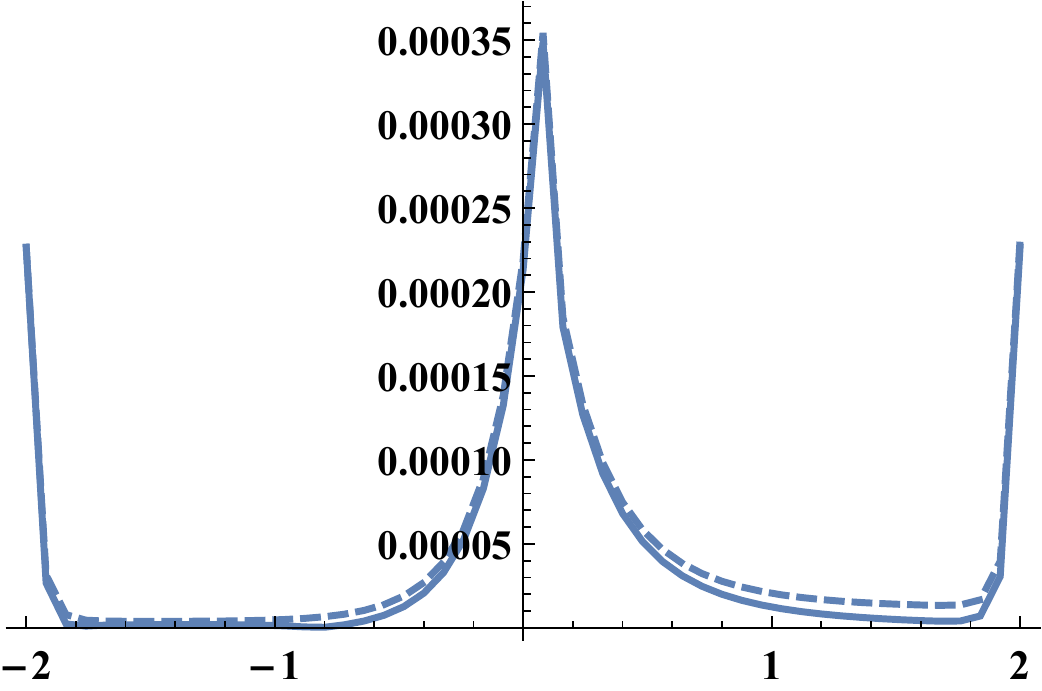}\label{fig2B}}
	\end{center}
	\begin{center}
		\caption{Error between exact and numerical solution (with the use of $50\times 50$ discretization points and after $5$ iterations steps).}\label{Fig-ErrBBL}
	\end{center}
\end{figure}

In Figure \ref{Fig-ErrBBL} the difference between exact solution of the stationary problem and the numerical solution (with the use of $50\times 50$ discretization points and after $5$ iterations steps) of the problem \eqref{main_num} at times $T=0.5$, $T=1$, $T=3$ are shown.

\end{example}

\section*{Acknowledgments}
This publication was made possible by NPRP grant NPRP 5-088-1-021
from the Qatar National Research Fund (a member of Qatar Foundation).
The statements made herein are solely the responsibility of the authors.

\section*{Conflict of Interests}
The authors declare that there is no conflict of interests regarding the publication of this paper.

\bibliographystyle{ieeetr}   
\bibliography{bubble}
\end{document}